\documentclass[11pt, reqno]{amsart}

 \usepackage{amsmath}
 \usepackage{amssymb}
\usepackage{bm}

\DeclareMathAccent{\mathring}{\mathalpha}{operators}{"17}

 \usepackage{color}

\newcommand{\mysection}[1]{\section{#1}
      \setcounter{equation}{0}}

\newtheorem{theorem}{Theorem}[section]
\newtheorem{lemma}[theorem]{Lemma}

\newtheorem{corollary}[theorem]{Corollary} 

\theoremstyle{definition}
\newtheorem{assumption}{Assumption}[section]

\theoremstyle{remark}
\newtheorem{remark}{Remark}[section]

\newtheorem{example}{Example}[section]

\newcommand{\tr}{\text{\rm tr}\,}
\newcommand{\loc}{\text{\rm loc}}

 \makeatletter 
 \def\dashint{%
 \operatorname%
 {\,\,\text{\bf--}\kern-.98em\DOTSI\intop\ilimits@\!\!}}
\def\ninf{\qopname\relax\@empty{inf\phantom{p}\!\!\!}}
 \makeatother

\newcommand\bbeta{\text{\raise-.2ex\hbox{$\bm{\beta}$}}}

\newcommand\esssup{\operatornamewithlimits{esssup}}

\newcommand\bR{\mathbb{R}}

\newcommand\bQ{\mathbb{Q}}
\newcommand\bS{\mathbb{S}}

\newcommand\cF{\mathcal{F}}

\newcommand\cN{\mathcal{N}}

\newcommand\frM{\mathfrak{M}}

\def\sft{{\sf t}}

\begin{document}

\title[Stochastic It\^o equations with drift in $L_{d+1}$]
{Some properties of
solutions of  It\^o equations 
with drift in $L_{d+1} $}

\author{N.V. Krylov}
 
\email{nkrylov@umn.edu}
\address{127 Vincent Hall, University of Minnesota,
 Minneapolis, MN, 55455}
 


\begin{abstract}
This paper is a natural continuation of \cite{Kr_20_2},
where strong Markov processes are constructed
in time inhomogeneous setting with Borel measurable
uniformly bounded and uniformly nondegenerate diffusion
and drift in $L_{d+1}(\mathbb{R}^{d+1})$. Here we study
some properties of these processes such as
 higher summability of 
Green's functions, boundedness of resolvent
operators in Lebesgue spaces, establish
 It\^o's formula, and so on.
\end{abstract}

\maketitle

\mysection{Introduction}

Let $\bR^{d}$ be a  Euclidean space of points
$x=(x^{1},...,x^{d})$, $d\geq 2$. Fix some 
 $p_{0},q_{0} \in[1,\infty)$  
such that
\begin{equation}
                                                    \label{5.10.1}
  \frac{d}{p_{0}}+\frac{1}{q_{0}}= 1.
\end{equation}

It is proved in \cite{Kr_20_2} that
It\^o's stochastic equations of the form
\begin{equation}
                                                 \label{11.29.2}
x _{t}=x  +\int_{0}^{t}\sigma (t +s,x_{s})\,dw_{s}
+\int_{0}^{t}b (t +s,x_{s}) \,ds
\end{equation}
admit  weak solutions,
where $w_{s}$ is a $d$-dimensional Wiener process,
$\sigma$ is a uniformly nondegenerate, bounded,
Borel function with values in the set of symmetric
$d\times d$ matrices, $b$ is a Borel measurable $\bR^{d}$-
valued function given on $\bR^{d+1}=(-\infty,\infty)\times\bR^{d}$
such that
\begin{equation}
                                                    \label{5.10.2}
\int_{\bR}\Big(\int_{\bR^{d}}|b(t,x)|^{p_{0}}\,dx\Big)^{q_{0}/p_{0}}\,dt<\infty 
\end{equation}
if $p_{0}\geq q_{0}$ or
\begin{equation}
                                                    \label{5.10.20}
\int_ {\bR^{d}} \Big(\int_{\bR}|b(t,x)|^{q_{0}}\,dt\Big)^{p_{0}/q_{0}}\,dx<\infty 
\end{equation}
if $ p_{0}\leq q_{0}$. 
Observe that the case $p_{0}=q_{0}=d+1$ is not excluded and in this case
the condition becomes $b\in L_{d+1}(\bR^{d+1})$.

The goal of this article is to study the properties of such
solutions or Markov processes whose trajectories
are solutions of \eqref{11.29.2}. In particular,
in Section \ref{section 10.25.2} for more or less general
processes of diffusion type we derive
several estimates of  Aleksandrov type by using
  Lebesgue spaces with mixed norms
like (in case $t=0$, $x=0$ in \eqref{11.29.2})
\begin{equation}
                                               \label{10.26.2}
E\int_{0}^{\infty}e^{- t}
f(t,x_{t}) \,dt\leq 
N \|  f\|_{L_{p ,q }},
\end{equation}
provided that $d/p+1/q\leq1$.

We also
show that expected time when the process $(t,x_{t})$,
starting at $(0,0)$, exits from $[0,R^{2})\times\{x:|x|<R\}$
is comparable to $R^{2}$. This plays a crucial role
in Section \ref{section 10.26.1} where we show a significant 
improvement of the Aleksandrov estimates in the direction
of lowering the powers of summability of $f$ in \eqref{10.26.2}
to $d_{0}/p+1/q\leq1$ with $d_{0}<d$. Time homogeneous
versions of these estimates are also given.

In the same Section \ref{section 10.25.2} we give some estimates
of the distribution of the exit times from cylinders,
which are heavily used in the sequel. In particular, we prove
that, for any $0\leq s\leq t<\infty$,
$$
E\sup_{r\in[s,t]}|x_{r}-x_{s}|^{ n}
\leq N|t-s|^{ n/2},
$$
regardless of wether $|t-s|$ is small or large.
It is to be said that instead of \eqref{5.10.2}
or \eqref{5.10.20}, which are not invariant
under self-similar transformations, we impose
a slightly stronger assumption on $b$,
that is invariant.

As we mentioned above, in Section \ref{section 10.26.1}
we improve the results of Section \ref{section 10.25.2}
in what concerns the Aleksandrov estimates, which allows us
to prove It\^o's formula for $ W^{1,2}_{p,q}(Q)$-functions
if $d_{0}/p+1/q\leq1$.  

Section \ref{section 10.24.1} is devoted to the case when our process
is, actually, not just of diffusion type, but a Markov
(time-inhomogeneous) process, whose existence is shown in
\cite{Kr_20_2}. We prove that their resolvent operators
are bounded in $L_{p,q}$.  
 
In the final Section \ref{section 12.21.1}
we prove the maximum principle for
the {\em transport \/} equation
$\partial_{t}u+b^{i}D_{i}u=0$ with $b\in L_{p,q}$
and $u$ from a Sobolev class. Results of that kind
seem to be new even if $u$ is smooth.

It is worth mentioning that there is a vast 
literature about stochastic and PDE equations when
\eqref{5.10.1} is replaced with $d/p+2/q\leq1$.
This condition is much stronger than ours.
Still we refer the reader to the recent articles
\cite{Na_18},
\cite{BFGM_19},  \cite{XZ_20}  and the 
references therein
for the discussion of many powerful 
and exciting results obtained
 under this stronger condition.
There are also many papers when this condition is
considerably relaxed on the account of
imposing various regularity conditions
on $\sigma$ and $b$ and/or considering
random initial conditions with bounded density,
 see, for instance,
\cite{ZZ_19}, \cite{Zh_20} and the 
references therein. Restricting
the situation to the one when $\sigma$ and $b$
are independent of time allows one to
relax the above conditions significantly
further, see, for instance, \cite{KS_19}
and the 
references therein.

 Introduce
$$
B_{R} =\{x\in\bR^{d}:| x|<R\}, 	\quad B_{R}(x)=x+B_{R},
\quad C_{T,R}=[0,T)\times B_{R},
$$
$$
C_{T,R}(t,x)=(t,x)+C_{T,R},\quad C_{R}(t,x)=C_{R^{2},R}(t,x),
\quad C_{R} =C_{R}(0,0),
$$
$$
 D_{i}=\frac{\partial}{\partial x^{i}},
\quad D_{ij}=D_{i}D_{j}\quad \partial_{t}=\frac{\partial}{\partial t}.
$$
For $p,q\in[1,\infty]$ and domains $Q\subset \bR^{d+1}$
 we introduce the space $L_{p,q}(Q)$ as the space of Borel
functions on $Q$ such that
$$
\|f\|^{q}_{L_{p,q}(Q)}:=
\int_{\bR}\Big(\int_{\bR^{d}}I_{Q}(t,x)|f(t,x)|^{p}\,dx\Big)^{q/p}\,dt<\infty 
$$
if $p\geq q$ or
$$
\|f\|^{p}_{L_{p,q}(Q)}:=\int_ {\bR^{d}} \Big(\int_{\bR}I_{Q}(t,x)
|f(t,x)|^{q}\,dt\Big)^{p/q}\,dx<\infty 
$$
if $ p\leq q$ with natural interpretation
of these definitions if $p=\infty$ or $q=\infty$.
If $Q=\bR^{d+1}$, we drop $Q$ in the above notation.
Observe that
$p $ is associated with   $x$ and
$q$ with   $t$ and the interior
integral is always elevated to the power $\leq 1$.
 In case $p=q=d+1$ we abbreviate $L_{d+1,d+1}=L_{d+1}$.
For the set of functions on $\bR^{d}$ summable to the $p$th
power we use the notation $L_{p}(\bR^{d})$.

If $\Gamma$ is a measurable subset of $\bR^{d+1}$ we denote by
$|\Gamma|$ its Lebesgue measure. The same notation
is used for measurable  subsets of $\bR^{d}$ with $d$-dimensional
Lebesgue measure. We hope that it will be clear
from the context which Lebesgue measure is used.
If $\Gamma$ is a measurable subset of $\bR^{d+1}$ and
$f$ is a function on $\Gamma$ we denote 
$$
\dashint_{\Gamma}f\,dxdt=\frac{1}{|\Gamma|}
\int_{\Gamma}f\,dxdt.
$$
In case $f$ is a function on a measurable subset $\Gamma$
of $\bR^{d}$ we set
$$
\dashint_{\Gamma}f\,dx =\frac{1}{|\Gamma|}
\int_{\Gamma}f\,dx .
$$

\mysection{The case of general diffusion type processes
with drift in $L_{p_{0},q_{0}}$}

                                        \label{section 10.25.2}

Let $(\Omega,\cF,P)$ be a complete probability
space, let $\cF_{t}, t\geq0$, be an increasing family of
complete $\sigma$-fields $\cF_{t}\subset\cF$,  
let $m_{t}$ be an $\bR^{d}$-valued continuous local martingale
relative to $\cF_{t}$, let  $A_{t}$ be a continuous
$\cF_{t}$-adapted
nondecreasing process, let $B_{t}$ be a continuous $\bR^{d}$-valued
$\cF_{t}$-adapted process which has finite variation (a.s.)
on each finite time interval and such that $d|B_{t}|\ll dA_{t}$.
 Assume that
$$
A_{0}=0,\quad m_{0}=B_{0}=0,\quad d\langle m\rangle_{t}\ll  dA_{t}
$$
and that we are also given progressively measurable relative to
$\cF_{t}$ nonnegative processes $r_{t}$ and $c_{t}$. Finally, 
introduce
$$
x_{t}= m_{t}+B_{t},\quad \sft_{t}=\int_{0}^{t}r_{s}\,dA_{s},
 \quad \phi_{t}=\int_{0}^{t}c_{s}\,dA_{s},
$$
$$
a^{ij}_{t}=\frac{1}{2}\frac{d\langle m^{i},m^{j}\rangle_{t}}{dA_{t}},
\quad b_{t}=\frac{dB_{t}}{dA_{t}}.
$$
 Also take a stopping time $\gamma$ and denote
$$
A=E\int_{0}^{\gamma}e^{-\phi_{t}} \tr a_{s}\,dA_{t},\quad 
B=E\int_{0}^{\gamma} e^{-\phi_{t}}\,|dB_{t}|.
$$

\begin{assumption}
                                         \label{assumption 10.20.1}
We are given a function $h\in L_{p_{0},q_{0},\loc}$
such that
$$
|b_{t}|\leq r^{1/q_{0}}_{t}(\det
a_{t} )^{1/p_{0}}h(\sft_{t}
,x_{t}).
$$
Furthermore, there exists a constant $\bar H>0$
and a continuous function $W(\theta)\geq0$, $\theta\in[0,1]$, such that
$W(0)=0  , W(1)=1$ and for any $(t,x)\in\bR^{d+1}$,
  $\rho>0$, and $\tau\in(0,\rho^{2}]$ we have
\begin{equation}
                                               \label{8.19.1}
 \|h\|^{q_{0}}_{L_{p_{0},q_{0}}(C_{\tau,\rho}(t,x))} \leq
\bar H \rho W(\tau/\rho^{2}).
\end{equation} 
\end{assumption}

\begin{remark}
                              \label{remark 12.20.1}
The reader will see that the assumption
that $W(0)=0$ can be replaced by the requirement that
$W(0)$ be sufficiently small, the smallness
defined only by $d,p_{0},\bar H, W$, and $\delta$,
where $\delta$ is introduced in Assumption
\ref{assumption 8.26.1} below.
\end{remark}

Observe that if $p_{0}=q_{0}=d+1$ and $h\in L_{d+2}$
(a typical case in the theory of parabolic equations),
then \eqref{8.19.1} is satisfied with $\bar H=
\|h\|^{d+1}_{L_{d+2}}$, because by H\"older's inequality
$$
\|h\|^{d+1}_{L_{d+1}(C_{\tau,\rho}(t,x))}
\leq \bar H\rho(\tau/\rho)^{2/(d+2)}.    
$$
 On the other hand, it may happen that
\eqref{8.19.1} is satisfied with $p_{0}=q_{0}=d+1$ 
but $h\not\in L_{d+2,\loc}$.

It is worth noting that a condition very similar to
\eqref{8.19.1} first appeared in \cite{Na_12}.  
 
\begin{example}
                                          \label{example 8.24.1}
Take  $\alpha\in(0,d),\beta\in(0,1)$ such that
$\alpha+2\beta=d+1$ and consider the function
$g(t,x)=|t|^{-\beta}|x|^{-\alpha}$. Observe that
$$
\int_{C_{\tau,\rho}(t,x)}g(s,y)\,dyds  
=\rho (\tau/\rho^{2})^{1-\beta}\int_{C_{1}(t',x')}g(s,y)\,dyds,
$$
where $t'=t/\tau$, $x'=x/\rho$. Obviously, the last integral
is a bounded function of $(t',x')$. Hence, the function
$h=g^{1/(d+1)}$ satisfies \eqref{8.19.1}
with $p_{0}=q_{0}=d+1 $. As is easy to see
for any $p>d+1$ one can find $\alpha$ and $\beta$
above such that $h\not\in L_{p,\loc}$.
 
\end{example}

The following is a particular case 
of Theorem 4.5 of \cite{Kr_20_2}. In Theorem \ref{theorem 9.27.1}
$p_{0}=\infty$ is allowed (and then $q_{0}=1$).
\begin{theorem}
                                        \label{theorem 9.27.1}

Suppose that  Assumption \ref{assumption 10.20.1}
is satisfied and
\begin{equation}
                                                  \label{10.4.3}
p,q\in[1,\infty],\quad\frac{d}{p }+\frac{1}{q }=1 .
\end{equation}
Then for any Borel $f\geq0$
\begin{equation}
                                                   \label{5.6.4}
  E\int_{0}^{\gamma}e^{-\phi_{t}} \kappa _{t}
f( \sft _{t},x_{t})\,dA_{t}\leq N(d,p_{0}  )
\Big(A+ \|h\|_{L_{p_{0},q_{0}}}^{2q_{0} }
\Big)^{d/(2p )}
\|f\|_{L_{p ,q }},
\end{equation}
where $\kappa _{t}=r_{t}^{1/q }(\det a_{t})^{1/p }$.

\end{theorem}

Our first goal is to estimate $A$ and eliminate it from 
\eqref{5.6.4}.
 
\begin{assumption}
                                   \label{assumption 8.26.1}
We have $\sft_{\infty}=\infty$ and
there is a constant $\delta\in(0,1)$
such that for all values of the arguments  
$$
\delta \,\tr a_{t}
\leq r _{t}.
 $$

\end{assumption}

This assumption 
and Assumption \ref{assumption 10.20.1}
are supposed to hold below in this section.
 
\begin{lemma}
                                      \label{lemma 8.16.1}
 
Introduce
$$
\tau_{R}=\inf\{t\geq 0:(\sft_{t},x_{t})\not\in C_{R}\} .
$$
Then
\begin{equation}
                                          \label{8.16.1}
A:=
E\int_{0}^{\tau_{R}}\tr a_{t}\,dA_{t}\leq \delta^{-1}R^{2},
\end{equation}
 and consequently (for $c_{t}\equiv 0$), assuming \eqref{10.4.3},
 for any Borel nonnegative $f$
\begin{equation}
                                          \label{9.29.2}
 E\int_{0}^{\tau_{R}} \kappa _{t}
f( \sft _{t},x_{t})\,dA_{t}\leq N(d,p_{0},  \delta)
 (1+\bar H^{d/p})R 
 ^{d/  p  }
\|f\|_{L_{p ,q }}.
\end{equation}
In particular (as $p=q=d+1$),
\begin{equation}
                                          \label{10.21.1}
 E\int_{0}^{\tau_{R}} (r_{t}\det a_{t})^{1/(d+1)}
 \,dA_{t}\leq N(d,p_{0}, \bar H,\delta)
 R  ^{2}.
\end{equation}
\end{lemma}

Proof. We have
$$
\int_{0}^{\tau_{R}}\tr a_{t}\,dA_{t}\leq \delta^{-1}
\int_{0}^{\tau_{R}}r_{t}\,dA_{t}=\delta^{-1}\sft_{\tau_{R}}
\leq \delta^{-1}R^{2} 
$$
and \eqref{8.16.1} follows.  After that \eqref{9.29.2}
follows from \eqref{8.19.1} and \eqref{5.6.4} with $c=0$,
 $f=I_{C_{R}}$, and $h$ replaced by $hI_{C_{R}}$.
The lemma is proved.

\begin{remark}
The detailed form of the right-hand side in
\eqref{9.29.2} will be important in 
Section \ref{section 12.21.1}.
\end{remark}

Estimate \eqref{10.21.1} says that in the typical case
of nondegenerate diffusion $\tau_{R}$ is of order
not more than $R^{2}$.
A very important fact which is implied by Corollary
\ref{corollary 7.29.1} below is that 
$\tau_{R}$ is of order
not less than $R^{2}$.
On the way of showing this we start with the following.

\begin{theorem}
                                      \label{theorem 8.2.1}
There is 
$\theta=\theta(d,\delta,p_{0},\bar H,W) >0$ such that for
all $R>0$
\begin{equation}
                                          \label{8.2.2}
P(\sft_{\tau_{R}}>  \theta R^{2} )\geq 1/2.
\end{equation}
\end{theorem}

Proof.  First we make an important claim, having a general character,
that
it suffices to prove \eqref{8.2.2} for $R=1$.
To prove the claim, fix an $R>0$ and  
make a self-similar transformation by setting
  $\hat x _{t}=x_{t}/R$,
$\hat \sft _{t}=\sft_{t}/R^{2}$.  
Then
$$
\hat x _{t}=\hat m _{t}+\hat B _{t},\quad  \hat \sft _{t}=\int_{0}^{t}
\hat r _{s}\,dA_{s},
$$
where
$$
\hat m _{t}=m_{t}/R,\quad \hat B _{t}=B_{t}/R,\quad
  \hat  r _{t}=r_{t}/R^{2}.
$$
Observe that
$$
\hat a^{ij}_{t}=\frac{1}{2}\frac{d\langle \hat m^{i},
\hat m^{j}\rangle_{t}}{dA_{t}}=a^{ij}_{t}/R^{2},\quad
\delta \,\tr\hat  a_{t}
\leq \hat r _{t}
$$
and for $\hat b_{t}=d\hat B_{t}/dA_{t}$ and $\hat h(t,x)
=Rh(R^{2}t,Rx)$
$$
|\hat b_{t}|=|b_{t}|/R\leq \hat r_{t}^{1/q_{0}}(\det \hat a_{t})
^{1/p_{0}}\hat h(\hat\sft_{t},\hat x_{t}).
$$
Furthermore, if $p_{0}\geq q_{0}$, for any $\rho>0$,
$\tau\leq \rho^{2}$, and $(t_{0},x_{0})\in\bR^{d+1}$    
$$
\int_{0}^{\tau}\Big(\int_{|x|\le \rho}  
|\hat h(t_{0}+ t,x_{0}+ x)|^{p_{0}}\,dx\Big)^{q_{0}/p_{0}}\,dt
$$
$$
=R^{q_{0}-2-q_{0}d/p_{0}}
\|h\|^{q_{0}}_{L_{p_{0},q_{0}}(C_{ R^{2}\tau,R\rho}(R^{2}t_{0},Rx_{0}) )}
$$
$$
\leq R^{q_{0}-2-q_{0}d/p_{0}}\bar H    R\rho W(\tau/\rho^{2}) 
=\bar H  \rho W(\tau/\rho^{2}).
$$

Also, if $q_{0}>p_{0}$,
$$
\int_{|x|\le \rho} \Big( \int_{0}^{\tau}
|\hat h(t_{0}+ t,x_{0}+ x)|^{q_{0}}\,dt\Big)^{p_{0}/q_{0}}\,dx
$$
$$
=R^{p_{0} -2p_{0}/q_{0}-d}
\|h\|^{p_{0}}_{L_{p_{0},q_{0}}(C_{ R^{2}\tau,R\rho}(R^{2}t_{0},Rx_{0}) )}
$$
$$
\leq R^{p_{0} -2p_{0}/q_{0}-d}\big(\bar H    R\rho W(\tau/\rho^{2})
\big)^{p_{0}/q_{0}} 
=\big(\bar H \rho W(\tau/\rho^{2})\big)^{p_{0}/q_{0}} .
$$

Hence, if \eqref{8.2.2} for $R=1$ is true, then
\begin{equation}
                                          \label{8.17.1}
P(\hat\sft_{\hat\tau_{1}}>  \theta )\geq 1/2,
\end{equation}
where
$$
\hat \tau_{1}=\inf\{t\geq0:(\hat \sft_{t},\hat x_{t})\not\in C_{1}
\}=\inf\{t\geq0:(  \sft_{t}, x_{t})\not\in C_{R}\}=\tau_{R},
$$
so that  \eqref{8.17.1} is just another form of \eqref{8.2.2}.

Having proved the claim we set $R=1$ and observe that by 
 \eqref{9.29.2}
$$
E\int_{0}^{\tau_{1}}|b(\sft_{t},x_{t})|I_{\sft_{t}\leq
\theta}\,dA_{t}\leq N(d,p_{0},  \delta)
 (1+\bar H^{d/p_{0}})
\|h\|_{L_{p_{0},q_{0}}(C_{\theta,1})}
$$
$$
\leq N(d,p_{0},  \delta)
 (1+\bar H^{d/p_{0}}) [\bar H W(\theta)]^{1/q_{0}} . 
$$
Furthermore, for $\gamma=\inf\{t\geq0:\sft_{t}\geq\theta\}$
$$
E\sup_{t\leq \gamma}|m_{t}|\leq 3
E\langle m\rangle_{\gamma}^{1/2}
\leq 6E\Big(\int_{0}^{\gamma}\tr a_{t}\,dA_{t}\Big)^{1/2}
$$
$$
\leq(6/\delta)E\Big(\int_{0}^{\gamma}r_{t}\,dA_{t}\Big)^{1/2}
\leq(6/\delta)\sqrt\theta.
$$

It follows that for $\theta$ small enough
$$
E\sup_{t\leq \tau_{1}\wedge\gamma}|x_{t}|\leq 1/2,
\quad P(\sup_{t\leq \tau_{1}\wedge\gamma}|x_{t}|\geq1)\leq 1/2,
\quad P(\sup_{t\leq \tau_{1}\wedge\gamma}|x_{t}|<1)\geq 1/2,
$$
and \eqref{8.2.2} follows. The theorem is proved.

\begin{corollary}
                                          \label{corollary 8.18.1}
For  $\lambda \in[0,1]$ and $R>0$ we have
\begin{equation}
                                          \label{8.18.3}
Ee^{-\lambda R^{-2}\sft_{\tau_{R}}}\leq e^{- \lambda \theta/4}.
\end{equation}
\end{corollary}

Indeed, self-similarity allows us to concentrate on $R=1$.
In that case
the derivative with respect to $\lambda$ of the left-hand
side of \eqref{8.18.3} is
$$
- E\sft_{\tau_{R}}e^{-\lambda\sft_{\tau_{1}}}\leq -
 e^{-\lambda  }
E\sft_{\tau_{1}}\leq-e^{-\lambda  }\theta/2,
$$
where the last inequality follows from \eqref{8.2.2}.
By integrating we find
$$
Ee^{-\lambda\sft_{\tau_{1}}}-1\leq  
(e^{-\lambda  }-1)\theta/2,
$$
which after using
$$
e^{-\lambda  }-1\leq- \lambda/2,\quad
1-\lambda\theta/4\leq e^{-\lambda\theta/4}
$$
leads to \eqref{8.18.3}.

\begin{theorem}
                                              \label{theorem 8.20.1}
For any $\lambda, R>0$ we have
\begin{equation}
                                          \label{8.20.1}
Ee^{-\lambda \sft_{\tau_{R}}}\leq e^{\theta/4}
 e^{- \sqrt\lambda R \theta/4}.
\end{equation}
In particular, for  
  any $R,t\in(0,\infty)$ we have
\begin{equation}
                                             \label{10.2.2}
P( \sft_{\tau_{R} }\leq  t )\leq 
e^{\theta/4} \exp\Big(-\frac{\theta^{2}R^{2}}{64 t}\Big).
\end{equation}

\end{theorem}

Proof. We may assume that $R=1$. In that case
take an integer $n\geq 1$, introduce $\tau^{k}$, $k=1,...,n$,
as the first exit time of $(\sft_{t},x_{t})$
from $C_{1/n}(\sft_{\tau^{k-1}},x_{\tau^{k-1}})$ after $\tau^{k-1}$
($\tau^{0}:=0$). In light of \eqref{8.18.3} we have
$$
E\Big(e^{-n^{2}(\sft_{\tau^{k}}-\sft_{\tau^{k-1} })}\mid \cF_{\tau^{k-1}}\Big)
\leq e^{-  \theta/4}.
$$
Hence
$$
Ee^{-n^{2} \sft_{\tau_{1}}}\leq E\prod_{k=1}^{n}
e^{-n^{2}(\sft_{\tau^{k}}-\sft_{\tau^{k-1} })}\leq  
e^{- n \theta/4},
$$
which for $n=\lfloor\sqrt{\lambda}\rfloor$
 yields \eqref{8.20.1}.
 
We prove \eqref{10.2.2} again assuming that $R=1$. We have
$$
P(\sft_{\tau_{1} }\leq  t  )=P\big(
 \exp(-\lambda 
\sft_{\tau_{1} }\geq \exp(-\lambda t)\big)
\leq e^{\theta/4} \exp( \lambda t-\sqrt{\lambda}\theta/4).
$$
For $\sqrt{\lambda}=\theta/(8t)$ we get \eqref{10.2.2}.
The theorem is proved.

\begin{corollary}
                                       \label{corollary 7.29.1}
There is a constant $N=N(\theta )$
such that for any $R\in(0,\infty)$ 
\begin{equation}
                                                   \label{8.21.1}
NE\sft_{\tau_{R}}\geq R^{2},\quad NE\int_{0}^{\tau_{R}}
e^{-\sft_{t}}r_{t}\,dA_{t}\geq R^{2}\wedge 1.
\end{equation}

\end{corollary}

Indeed, for any $t_{0}>0$
$$
E\sft_{\tau_{R}}\geq  t_{0}R^{2}P(\sft_{\tau_{R}}>t_{0}R^{2})\geq 
R^{2}
t_{0}\Big(1-e^{\theta/4} \exp\Big(-\frac{\theta^{2} }{64 t_{0}}\Big)\Big),
$$
$$
E\int_{0}^{\tau_{R}}
e^{-\sft_{t}}r_{t}\,dA_{t}=E(1-e^{-\sft_{\tau_{R}}})
\geq EI_{\sft_{\tau_{R}} >t_{0}R^{2}} (1-e^{-t_{0}R^{2}})
$$
$$
=P(\sft_{\tau_{R}}>t_{0}R^{2})(1-e^{-t_{0}R^{2}})
$$
$$
\geq\Big(1-e^{\theta/4} \exp\Big(-\frac{\theta^{2} }{64 t_{0}}\Big)\Big)
(1-e^{-t_{0}R^{2}}),
$$
which yields \eqref{8.21.1} for an appropriate $t_{0}$.

\begin{corollary}
                                        \label{corollary 10.26.1}
Assume that $\sft_{t}\leq \delta^{-1}t$
for all $t\geq 0$. Then for any $n>0$ and
  $0\leq s\leq t $ we have
\begin{equation}
                                            \label{10.28.2}
E\sup_{r\in[s,t]}|x_{r}-x_{s}|^{ n}
\leq N|t-s|^{ n/2},
\end{equation}
where $N=N(n,\delta,\theta)$.
\end{corollary}

Indeed, clearly we may assume that $s=0$. Then
$$
E\sup_{r\leq t]}|x_{r} |^{ n}
=n\int_{0}^{\infty}\mu^{n-1}
P(\sup_{r\leq t }|x_{r} |\geq \mu)\,d\mu
\leq n\int_{0}^{t^{1/2}}\mu^{n-1}\,d\mu
$$
$$
+e^{\theta/4}\int_{t^{1/2}}^{\infty}\mu^{n-1}
\exp\Big(-\frac{\delta\theta^{2}\mu^{2}}{64 t}\Big)\,d\mu,
$$
where the inequality is due to
\eqref{10.2.2} and the fact that, for $\mu> t^{1/2}$,
$$
\{\sup_{r\leq t }|x_{r} |\geq \mu\}\subset\{\tau_{\mu}\leq t\}
\subset\{\sft_{\tau_{\mu}}\leq t/\delta\}.
$$
Changing variables $\mu=\nu\sqrt t$ leads to \eqref{10.28.2}.

\begin{remark}
Observe that the right-hand side of \eqref{10.28.2}
has the same form for small $|t-s|$ and for large ones.
\end{remark}

Now we start moving toward estimating the resolvent
of Markov diffusion processes in $L_{p,q}$.
\begin{lemma}
                                            \label{lemma 8.22.1}
Assume \eqref{10.4.3}.
Then there is a constant $N$,
depending only on
$\delta,d,\bar H,p, p_{0},W$, such that
for any   $t_{0}\geq0$, $x_{0}\in \bR^{d}$,
and Borel nonnegative $f$ vanishing outside $C_{1}(t_{0},x_{0})$
we have 
\begin{equation}
                                                   \label{8.22.10}
 E\int_{0}^{\infty}e^{- \sft_{t}}
f(\sft_{t},x_{t})\kappa _{t} \,dA_{t}\leq 
N\Phi(t_{0},x_{0})\|  f\|_{L_{p ,q }},
\end{equation}
where $\kappa _{t}=r_{t}^{1/q }(\det a_{t})^{1/p }$
and $\Phi (t,x)=\exp(-  (|x|+ \sqrt t)\theta/8)$.
\end{lemma}

Proof.  
Fix $\rho=\rho(\theta)>0$ such that the right-hand side
of \eqref{10.2.2} equals 1/2 when $R=\rho$ and $t=1$. Then
introduce
  $\tau^{0}$ as the first time $(\sft_{t},x_{t})$
hits $\bar C_{1}(t_{0},x_{0})$ and set $\gamma^{0}$
as the first time after $\tau^{0}$ the process
$(\sft_{t},x_{t})$ exits from $C_{1,1+\rho}(t_{0},x_{0})$.
We define recursively $\tau^{k}$, $k=1,2,...$, as the first
time after $\gamma^{k-1}$ the process
$(\sft_{t},x_{t})$
hits $\bar C_{1}(t_{0},x_{0})$ and $\gamma^{k}$ as 
the first time after $\tau^{k}$ the process
$(\sft_{t},x_{t})$ exits from $C_{1,1+\rho}(t_{0},x_{0})$.

 These stopping times are either
infinite or lie between $t_{0}$ and $t_{0}+1$.
Therefore  the left-hand side of
\eqref{8.22.10} equals
\begin{equation}
                                                   \label{8.22.1}
 E
\sum_{k=0}^{\infty} e^{- \sft_{\tau^{k}}}I_{k},
\end{equation}
where
$$
I_{k}=I_{\tau^{k}\in[t_{0},t_{0}+1]}
E\Big(\int_{\tau^{k}\wedge(t_{0}+1)}^{ \gamma^{k}\wedge(t_{0}+1)}
e^{-( \sft_{t}- \sft_{\tau^{k}})}
f(\sft_{t},x_{t})\kappa _{t}\,dA_{t} \mid \cF_{\tau^{k}}\Big).
$$
Here  
$$
\int_{\tau^{k}\wedge(t_{0}+1)}^{ \gamma^{k}\wedge(t_{0}+1)}
e^{-( \sft_{t}- \sft_{\tau^{k}})}
\tr a_{t}\,dA_{t}\leq\delta^{-1}
\int_{\tau^{k}\wedge(t_{0}+1)}^{ \gamma^{k}\wedge(t_{0}+1)}
 r_{t}\,dA_{t}
$$
$$
=\delta^{-1}(\sft_{\gamma^{k}\wedge(t_{0}+1)}
-\sft_{\tau^{k}\wedge(t_{0}+1)})\leq \delta^{-1}.
$$
It follows  from \eqref{5.6.4}
that
$I_{k}\leq N\|f\|_{L_{  p , q }}$.  

Next, observe that, if $\sqrt{t_{0}}>|x_{0}| $, then  $\tau^{0}$
is bigger than the first exit time of $(\sft_{t},x_{t})$
from $C_{\sqrt{t_{0}}}$, and
by Theorem \ref{theorem 8.20.1} 
$$
Ee^{- \sft_{\tau^{0}}}\leq Ne^{-\sqrt{t_{0}}\theta/4}.
$$
In case $\sqrt{t_{0}}\leq |x_{0}| $ and $|x_{0}|>1$
our $\tau^{0}$
is bigger than the first exit time of $(\sft_{t},x_{t})$
from $C_{|x_{0}|-1}$, and
$$
Ee^{- \sft_{\tau^{0}}}\leq Ne^{-(|x_{0}|-1)\theta/4}
\leq Ne^{-|x_{0}|\theta/4}.
$$
The last estimate (with different $N$) also holds if $|x_{0}|\leq1$
and we conclude that in all cases
$$
Ee^{- \sft_{\tau^{0}}}\leq Ne^{-(\sqrt{t_{0}}+|x_{0}|)\theta/8}.
$$

Furthermore,  by the choice of $\rho$
$$
E\Big(e^{- (\sft_{\gamma^{k}}
-\sft_{\tau^{k}})}\mid \cF_{\tau^{k}}\Big)\leq \frac{1}{2},
$$
$$
Ee^{- \sft_{\tau^{k}}}=
Ee^{- \sft_{\gamma^{k-1}}}E\Big(e^{- (\sft_{\tau^{k}}
-\sft_{\gamma^{k-1}})} \mid \cF_{\gamma^{k-1}}\Big)
\leq \frac{1}{2}Ee^{- \sft_{\gamma^{k-1}}},
$$
so that
$$
Ee^{- \sft_{\tau^{k}}}\leq\frac{1}{4}
Ee^{-  \sft_{\tau^{k-1}}},\quad
Ee^{-  \sft_{\tau^{k}}}\leq 4^{-k}Ee^{-  \sft_{\tau^{0}}}.
$$

Recalling \eqref{8.22.1} we see that the left-hand side 
of \eqref{8.22.10}
is dominated by
$N\Phi(t_{0},x_{0})\|f\|_{L_{p ,q }}$  and the lemma is proved.

 The following theorem shows that the time spent
by $(\sft_{t},x_{t})$ in cylinders $C_{1}(0,x)$ far from the origin
decays very fast.
\begin{theorem}
                                            \label{theorem 8.22.1}
Suppose that    
\begin{equation}
                                                   \label{10.21.3}
p ,q \in[1,\infty],\quad
\nu:=1-\frac{d}{p }-\frac{1}{q }\geq0 .
\end{equation}
Then
 there is a constant $N=N(\delta,d,\bar H,p,q,p_{0},W)$ such that
for any $\lambda>0$ and   Borel nonnegative $f$  
we have 
\begin{equation}
                                                   \label{8.22.4}
 E\int_{0}^{\infty}e^{- \lambda\sft_{t}}
f(\sft_{t},x_{t})\kappa _{t}r_{t}^{\nu}\,dA_{t}\leq 
N\lambda^{-\nu-d/(2p )}\|\Phi_{\lambda}^{1-\nu} f\|_{L_{p ,q }},
\end{equation}
where $\Phi _{\lambda}(t,x)=\exp(- \sqrt\lambda 
(|x|+ \sqrt t)\theta/32)$ and $r_{t}^{\nu}:=1$ if $\nu=0$.
\end{theorem}

Proof. First assume that $\nu=0$.
The self-similarity argument allows us to also assume that 
$\lambda=1$. In that case take a nonnegative $\zeta\in C^{\infty}_{0}
(\bR^{d+1})$ with support in $C_{1}$ and
 unit integral and for $(t,x),(s,y)\in\bR^{d+1}$
set
$$
f_{s,y}(t,x)=f(t,x)\zeta(t-s,x-y).
$$
Clearly, due to Lemma \ref{lemma 8.22.1},
$$
I:=E\int_{0}^{\infty}e^{-  \sft_{t}}
f(\sft_{t},x_{t})\kappa _{t}\,dA_{t}
=\int_{0}^{\infty}\int_{\bR^{d }}E\int_{0}^{\infty}e^{-  \sft_{t}}
f_{s,y}(\sft_{t},x_{t})\kappa _{t}\,dA_{t}\,dyds
$$
$$
\leq N\int_{0}^{\infty}\int_{\bR^{d }}
\Phi(s,y)\|  f_{s,y}\|_{L_{p ,q }}\,dyds.
$$
{\em Case $p\geq q$}. Then $q<\infty$ and
it follows by H\"older's 
inequality  that
$$
I\leq N\Big(
\int_{0}^{\infty}\int_{\bR^{d }}
\Phi^{q /2}(s,y)\int_{0}^{\infty}\Big(\int_{\bR^{d}}
f_{s,y}^{p }(t,x)\,dx\Big)^{q /p }dt
\,dyds\Big)^{1/q }
$$
$$
=N\Big(
\int_{0}^{\infty}dt \Big(\int_{0}^{\infty}
\int_{\bR^{d }}\Phi^{q /2}(s,y)
 \Big(\int_{\bR^{d}}
f_{s,y}^{p }(t,x)\,dx\Big)^{q /p } 
\,dyds\Big)\Big)^{1/q } 
$$

$$
\leq N\Big(
\int_{0}^{\infty}dt\Big(\int_{0}^{\infty}\int_{\bR^{d }}\int_{\bR^{d }}
\Phi^{p /4}(s,y)f_{s,y}^{p }(t,x)\,dydsdx 
\Big)^{q /p }\Big)^{1/q }.
$$

We replace $\Phi^{p /4}(s,y)$ by $\Phi^{p /4}(t,x)$
taking into account that these values are comparable
as long as $\zeta(t-s,x-y)\ne0$. After that
integrating over $dyds$ leads immediately to
\eqref{8.22.4}. 

{\em Case $p< q$}. 
It follows by H\"older's inequality that  
$$
I\leq N\Big(
\int_{0}^{\infty}\int_{\bR^{d }}
\Phi^{p /2}(s,y)\int_{\bR^{d}}\Big(\int_{0}^{\infty}
f_{s,y}^{q}(t,x)\,dt\Big)^{p/q  }dx
\,dyds\Big)^{1/p }
$$
$$
\leq N\Big(\int_{\bR^{d}}dx\Big(
\int_{0}^{\infty}\int_{\bR^{d }}\int_{0}^{\infty}
\Phi^{q /4}(s,y)f_{s,y}^{q}(t,x)\,dtdyds\Big)^{p/q}
\Big)^{1/p}.
$$
This leads   to
\eqref{8.22.4} as above.
The theorem is proved if $\nu=0$.

If $\nu\in(0,1)$,
by H\"older's inequality the left-hand side of
\eqref{8.22.4} is dominated by $I_{1}I_{2}$, where
$$
I_{1}^{1/\nu}
= E\int_{0}^{\infty}e^{-\lambda \sft_{t}}r_{t}\,dA_{t}
= 1/\lambda,
$$
$$
I_{2}^{1/(1-\nu)}=
E\int_{0}^{\infty}e^{-\lambda \sft_{t}}r^{1/(q-\nu q)}_{t}
(\det a_{t})^{1/(p-\nu p)}f^{1/(1-\nu)}(\sft_{t},x_{t})\,dA_{t}.
$$
Here
$$
\frac{d}{p-p\nu}+\frac{1}{q-q\nu}=1
$$
so that by the case that $\nu=0$
$$
I_{2}^{1/(1-\nu)}\leq N\lambda^{-d/(2p-2p\nu)}\|\Phi_{\lambda}
f^{1/(1-\nu)}\|_{L_{p-p\nu,q-q\nu}}
$$
$$
=N\lambda^{-d/(2p-2p\nu)}\|\Phi_{\lambda}^{1-\nu}f\|_{L_{p,q}}
^{1/(1-\nu)}.
$$
This leads   to
\eqref{8.22.4} again. Finally, if $\nu=1$ so that
$p=q=\infty$, the left-hand side of \eqref{8.22.4}
is obviously dominated by $\lambda^{-1}\sup f$, so that
\eqref{8.22.4} holds with $N=1$. The theorem is
proved.

By taking $q=\infty$ and $f(t,x)=f(x)$ we come to the following,
which extends Corollary 2.4 of \cite{Kr_19_1} to the case
of time dependent drift $b\in L_{p_{0},q_{0},\loc}$.
It is further generalized in a more specific situation
by relaxing the restriction on $p$
in Theorem \ref{theorem 8.30.1}.
\begin{corollary}
                                      \label{corollary 10.4.1}
Let $p \in [d,\infty]$.
Then   
for any $\lambda>0$ and   Borel nonnegative $f(x)$  
we have 
\begin{equation}
                                                   \label{10.4.2}
 E\int_{0}^{\infty}e^{- \lambda\sft_{t}}
f( x_{t})(\det a_{t})^{1/p}r_{t}^{1-d/p}\,dA_{t}\leq 
N\lambda^{- 1+d/(2p )}\|\Psi_{\lambda}^{d/p} f\|_{L_{p}(\bR^{d})},
\end{equation}
where $\Psi _{\lambda}( x)=\exp(- \sqrt\lambda  |x| \theta/32)$
and $N=N(\delta,d,\bar H,p  ,p_{0},W)$.
\end{corollary}  
 
Next results are dealing with the exit times of the process
$x_{t}$ rather than $(\sft_{t},x_{t})$. We will need them
while showing an improved integrability of Green's functions.
 
\begin{corollary}
                                      \label{corollary 9.29.1}
 
Suppose that  
\begin{equation}
                                          \label{9.29.4}
\delta r_{t}^{d}\leq\det a_{t}
\end{equation}
and for $|x|<R$ introduce
$$
  \tau'_{R}(x)=\inf\{t\geq 0: x+ x_{t} \not\in B_{R}\},
\quad \tau'_{R} =\tau'_{R}(0).
$$
Then there exists a constant $T=T(d,p_{0}, \bar H,\delta)$
such that
\begin{equation}
                                          \label{9.29.3}
P(\sft_{\tau'_{R}(x)}\geq TR^{2})\leq 1/2.
\end{equation}

\end{corollary}

To prove this, observe that $\tau'_{R}(x)\leq  
\tau'_{2R}$,
which shows that it suffices to prove 
\eqref{9.29.3} for $x=0$. In that case  
for $\varepsilon\in(0,1]$ introduce
$$
 \tau'_{\varepsilon,R}=\inf\{t\geq 0: 
(\sft_{t},x_{t}) \not\in C_{R^{2},\varepsilon R}\}.
$$
Then \eqref{9.29.2} with $f$ being the indicator
of $C_{R^{2},\varepsilon R}$ implies that
$$
E\sft_{\tau'_{\varepsilon,R}}\leq N
 E\int_{0}^{\tau_{R}} \kappa _{t}
f( \sft _{t},x_{t})\,dA_{t}\leq N\varepsilon^{d/p_{0}} R^{2}.
$$
For an appropriate choice of $\varepsilon$ the last
term is less than $(1/2)R^{2}$ and then
$$
P(\sft_{\tau'_{\varepsilon}}\geq R^{2})\leq
P(\sft_{\tau'_{\varepsilon,R}}\geq R^{2})\leq 1/2.
$$
Self-similarity  implies that \eqref{9.29.3}
holds with $T=\varepsilon^{-2}$ if $x=0$.

Estimate \eqref{10.4.7} below in case $b$ is bounded
was the starting point for the theory
of {\em time homogeneous\/} controlled diffusion processes.
 \begin{lemma}
                                      \label{lemma 9.29.1}
 Assume \eqref{9.29.4} and \eqref{10.4.3}.
Then
\begin{equation}
                                          \label{8.2.1}
E\int_{0}^{\tau'_{R}} 
 \tr a_{t}\,dA_{t}\leq \delta^{-1}E\sft_{\tau'_{R} }
\leq N(\delta,d,p_{0}, 
\bar H)R^{2},
\end{equation}
 and  for any Borel nonnegative $f(x)$ 
\begin{equation}
                                          \label{9.29.5}
 E\int_{0}^{\tau'_{R}}  \kappa _{t}
f( x_{t})\,dA_{t}\leq N(d,p_{0}, p,\bar H,\delta)
 R 
 ^{2-d/  p  } 
\|f\|_{L_{p  }(\bR^{d})}.
\end{equation}
In particular (as $p=d,q=\infty$ in $\kappa_{t}$),
\begin{equation}
                                          \label{10.4.7}
 E\int_{0}^{\tau'_{R}} (\det a_{t})^{1/d}
f( x_{t})\,dA_{t}\leq N(d,p_{0}, \bar H,\delta)
 R 
\|f\|_{L_{d}(\bR^{d})}.
\end{equation}
\end{lemma}

Proof. The first inequality in \eqref{8.2.1} is obvious.
To prove the second one observe that with $T$
from Corollary \ref{corollary 9.29.1}
$$
E\sft_{\tau'_{R} }\leq\sum_{k=1}^{\infty}
kTR^{2}P\big((k-1)TR^{2}\leq \sft_{\tau'_{R} }<kTR^{2})\big.
$$
For the event $\{(k-1)TR^{2}\leq \sft_{\tau'_{R} }<kTR^{2}\}$
to realize, the first coordinate 
of the process $(\sft_{t},x_{t})$ should
experience $k-1$ moves from $iTR^{2}$ to $(i+1)TR^{2}$
before $x_{t}$ exits from $B_{R}$. According to 
a ``conditional'' version of
Corollary \ref{corollary 9.29.1} the conditional
probability of each such move given the past
is less than $1/2$. Therefore,
$$
E\sft_{\tau'_{R} }\leq\sum_{k=1}^{\infty}
kTR^{2}2^{-k+1}
$$
and this yields the second inequality in \eqref{8.2.1}.
 
To prove \eqref{9.29.5} observe that for $k=1,2,...$,
$\Delta_{k} =[(k-1)TR^{2},kTR^{2})$,
$$
f_{k}(t,x)=I_{\Delta_{k}}(t)f(x)
$$ 
and $\gamma_{k}$ defined as the first time $\sft_{t}$
hits $kTR^{2}$,
according to
 \eqref{5.6.4} with $c=0$ on the set $ \tau'_{R} 
\geq \gamma_{k-1}$ we have
$$
E\Big(\int_{\gamma_{k-1}}^{\gamma_{k}
\wedge \tau'_{R}}\kappa_{t}f(x_{t})\,dA_{t}
\mid \cF_{\gamma_{k-1}}\Big)
=E\Big(\int_{\gamma_{k-1}}^{\gamma_{k}
\wedge \tau'_{R}}\kappa_{t}f_{k}(\sft_{t},x_{t})\,dA_{t}
\mid \cF_{\gamma_{k-1}}\Big)
$$
$$
\leq N\Big(R^{2}+\|h\|_{L_{p_{0},q_{0}}(\Delta_{k}\times B_{R})}^{2q_{0} }
\Big)^{d/(2p )}R^{2/q}
\|f\|_{L_{p  }(\bR^{d})}\leq NR^{2-d/p}\|f\|_{L_{p  }(\bR^{d})}.
$$
It follows that
$$
E \int_{0}^{ \tau'_{R}}\kappa_{t}f(x_{t})\,dA_{t}
=\sum_{k=1}^{\infty}
E \int_{\gamma_{k-1}\wedge \tau'_{R}}^{\gamma_{k}
\wedge \tau'_{R}}\kappa_{t}f(x_{t})\,dA_{t}
$$
$$
=\sum_{k=1}^{\infty}
E I_{\tau'_{R}\geq \gamma_{k-1}}\int_{\gamma_{k-1} }^{\gamma_{k}
\wedge \tau'_{R}}\kappa_{t}f(x_{t})\,dA_{t}
$$
$$
\leq NR^{2-d/p}\|f\|_{L_{p  }(\bR^{d})}\sum_{k=1}^{\infty}
P(\tau'_{R}\geq \gamma_{k-1}).
$$
As we have seen above, each of the probabilities in the last sum
is less than $2^{-k+1}$ and this proves the lemma.

 \mysection{Green's functions}
                                         \label{section 10.26.1}

In this section we make one step closer
to stochastic It\^o equations and restrict the class
of processes from Section \ref{section 10.25.2}.
We suppose that
$$
A_{t}=t,\quad  r_{t}\equiv 1,\quad \delta^{-1}|\lambda|^{2}\geq
(a_{t}\lambda,\lambda)\geq \delta|\lambda|^{2},
$$
Assumption  \ref{assumption 10.20.1} is satisfied
and also suppose that we are given 
$p$ and $q$ satisfying \eqref{10.21.3}.

Here is a straightforward consequence
of \eqref{8.22.4}.

\begin{theorem}
                                       \label{theorem 9.3.1}
 
There exists  a constant  $N=N(\delta,d,\bar H,p,q,p_{0},W)$ such that
for any $\lambda>0$ 
  there exists a nonnegative Borel function $G_{\lambda}(t,x)$
(Green's function of $(\sft_{\cdot},x_{\cdot})$)
on $\bR^{d+1}$
 such that $G_{\lambda}(t,x)=0$ for $t\leq0$
and for any Borel nonnegative $f$ given on $\bR^{d+1}$
we have
$$
 E\int_{0}^{\infty}e^{-\lambda t}
f(\sft_{t},x_{t}) \,d t =\int_{\bR^{d+1}}f(t,x)G_{\lambda}(t,x)\,dxdt,
$$
\begin{equation}
                                                   \label{9.3.5}
 \|\Phi^{\nu-1}_{\lambda} G_{\lambda}\|_{L'_{p ,q } } 
\le N\lambda^{-\nu-d/(2p )},
\end{equation}
where we use the notation
$$
\|u\|_{L'_{p,q} }=\Big(\int_{\bR}\Big(
\int_{\bR^{d}}|u(t.x)|^{p '}
\,dx\Big)^{q '/p '}\,dt\Big)^{1/q '}
\quad\text{if}\quad p\geq q,
$$
$$
\|u\|_{L'_{p,q} }=\Big(\int_{\bR^{d}}\Big(
\int_{\bR }|u(t.x)|^{q '}
\,dt\Big)^{p '/q '}\,dx\Big)^{1/p'}
\quad\text{if}\quad p< q,
$$
and $p '=p /(p -1),q '=q /(q -1)$.
 
\end{theorem}

The highest power of pure ($p =q =d+1$) 
summability of $G_{\lambda}$
 guaranteed by this theorem
 is $1+1/d$. It turns out that,
actually,
$G_{\lambda}$ is summable to a higher power. The proof of this
is based on the parabolic version of
Gehring's lemma from \cite{GS_82}.

Introduce $\bQ$ as the set of cylinders $C_{R}(t,x)$, $R>0$,  
$t\geq0$, $x\in\bR^{d}$. For $Q=C_{R}(t,x)\in\bQ$ let
 $2Q=C_{2R}(t,x)$. If $Q\in\bQ$ and $Q=C_{R}(t,x)$,
we call $R$ the radius of $Q$.

\begin{theorem}
                                       \label{theorem 9.3.2}
There exist  $d_{0}\in(1,d)$ and a constant    $N $, depending only
on $d,\delta$, $\bar H$, and $W$, such that for any     $Q\in\bQ$
of radius $R\leq 1/2$ and $p\geq d_{0}+1 $, we have
\begin{equation}
                                          \label{10.14.1}
\| G_{1}\|_{L_{p/(p-1)}(Q)}\leq N R^{-(d+2) /p }
\| G_{1}\|_{L_{1}(2Q )} , 
\end{equation}
which is equivalently rewritten as
$$
\Big(\dashint_{Q}G^{p/(p-1)}_{1}\,dxdt\Big)^{(p-1)/p}
\leq N\dashint_{2Q }G_{1}\,dxdt.
$$
 
\end{theorem}

Proof. We basically follow the idea in \cite{FS_84}.
Take $R\in(0,1/2]$,    $Q\in \bQ$ of radius $R$, and
define 
recursively 
$$
\gamma^{1}=\inf\{t\geq0 :(\sft_{t},x_{t})\in \bar Q\},\quad
\tau^{1}=\inf\{t\geq\gamma^{1} :(\sft_{t},x_{t})\not\in 2Q\},
$$
$$
\gamma^{n+1}=\inf\{t\geq\tau^{n} :(\sft_{t},x_{t})\in \bar Q\},\quad
\tau^{n+1}=\inf\{t\geq\gamma^{n+1} :(\sft_{t},x_{t})\not\in 2Q\}.
$$
Then for any nonnegative Borel $f$ vanishing outside $Q$
with $\|f\|_{L_{d+1}(Q)}=1$
we have
$$ 
\int_{Q}f   G_{1}(t,x)\,dxdt
=E\int_{0}^{\infty}e^{-t}f (t,x_{t})\,dt 
$$
$$
=\sum_{n=1}^{\infty}Ee^{- \gamma^{n}} 
E\Big(\int_{\gamma^{n}}^{\tau^{n}}e^{-(t-
 \gamma^{n} )}f(t,x_{t})\,dt
\mid \cF_{\gamma^{n}}\Big).
$$
Next we use the conditional  version of 
\eqref{9.29.2}
 to see that the conditional expectations
above are less than $NR^{d/(d+1)}$. After that
we use the conditional  version of Corollary
\ref{corollary 7.29.1}
to get that
$$
R^{2}\leq N E\Big(\int_{\gamma^{n}}^{\tau^{n}}e^{-(t-\gamma^{n})} \,dt
\mid \cF_{\gamma^{n}}\Big).
$$
Then we obtain
$$
\int_{Q}f   G_{1}(t,x)\,dxdt
\leq N  R ^{-(d+2)/(d+1)}
\sum_{n=1}^{\infty}E e^{-\gamma^{n}}
E\Big(\int_{\gamma^{n}}^{\tau^{n}}e^{-(t-\gamma^{n})} \,dt
\mid \cF_{\gamma^{n}}\Big)
$$
$$
=N  R ^{-(d+2)/(d+1)}
\sum_{n=1}^{\infty}E  
 \int_{\gamma^{n}}^{\tau^{n}}e^{- t } \,dt
$$
$$
\leq N  R ^{-(d+2)/(d+1)}
E \int_{0}^{\infty}e^{-t}I_{2Q}(t,x_{t}) \,dt
=N  R ^{-(d+2)/(d+1)}\int_{2Q}G_{1}(t,x)\,dxdt.
$$

The arbitrariness of $f$ 
implies that
$$
\Big(\dashint_{Q}  G_{1} ^{(d+1)/d}(t,x)\,dxdt
\Big)^{d/(d+1)}\leq N \dashint_{2Q}  G_{1}(t,x)\,dxdt.
$$

Now the assertion of the theorem for $p=d_{0}$ follows directly from
the corrected version of the famous Gehring's lemma
stated as Proposition 1.3  in \cite{GS_82}. 
For larger $p$ it suffices to use H\"older's inequality.
The theorem is proved.

\begin{remark}
The result of Theorem \ref{theorem 9.3.2}
in case $x_{t}$ is a solution of a {\em stochastic equation\/}
with deterministic $a(t,x)$ and $b\equiv0$
can be found in the PDE terms in Remark 1
of \cite{Es_93}, where it is given without proof.
Complete proof for bounded $b(t,x)$
can be extracted from \cite{CKKS_98} dealing with viscosity  
solutions of fully nonlinear parabolic equations.
With little effort this proof can be used 
to obtain the result of Theorem \ref{theorem 9.3.2}
if $b$ is bounded. Our $b$ is in $L_{p,q}$.
\end{remark}

Similar improvement of integrability occurs for
the Green's function of $x_{t}$ rather than $(t,x_{t})$.
Here is a straightforward consequence
of \eqref{10.4.2}.  

\begin{theorem}
                                       \label{theorem 10.4.1}
 Let $p\in[d,\infty)$. Then
for any $\lambda>0$ 
  there exists a nonnegative Borel function $g_{\lambda}( x)$
(Green's function of $ x_{\cdot} $)
on $\bR^{d }$
 such that   for any Borel nonnegative $f$ given on $\bR^{d }$
we have
$$
 E\int_{0}^{\infty}e^{-\lambda t}
f( x_{t}) \,d t =\int_{\bR^{d }}f( x)g_{\lambda}( x)\,dx ,
$$
\begin{equation}
                                                   \label{10.4.5}
 \|\Psi^{-d/p}_{\lambda} g_{\lambda}\|_{L_{p' } } 
\le N\lambda^{-1+d/(2p )},
\end{equation}
where   $\Psi _{\lambda}( x)=\exp(- \sqrt\lambda  |x| \theta/32)$,
  $p '=p /(p -1)$, and   $N$
depends only on $\delta$, $d$, $\bar H,p, p_{0},W$.
 
\end{theorem}

According to this theorem 
this Green's function is summable to the power $d/(d-1)$.
Again it turns out that this power can be increased.
If $B$ is an open ball in $\bR^{d}$ by $2B$ we denote
the concentric open ball of twice the radius of $B$.

\begin{theorem}
                                           \label{theorem 10.4.3}
There exist  $d_{0}\in(1,d)$ and a constant    $N $, depending only
on $d,\delta$, $\bar H$, and $W$, such that for any   ball  $B$
of radius $R\leq 1/2$ and $p\geq d_{0}  $, we have
\begin{equation}
                                          \label{10.4.6}
\| g_{1}\|_{L_{p/(p-1)}(B)}\leq N R^{-d /p }
\| g_{1}\|_{L_{1}(2B )} , 
\end{equation}
which is equivalently rewritten as
$$
\Big(\dashint_{B}g^{p/(p-1)}_{1}\,dx \Big)^{(p-1)/p}
\leq N\dashint_{2B }g_{1}\,dx.
$$
 
\end{theorem}

Proof. We again follow the idea in \cite{FS_84}.
Take $R\in(0,1/2]$, an open ball  $B$ of radius $R$, and
define 
recursively 
$$
\gamma^{1}=\inf\{t\geq0 : x_{t} \in \bar B\},\quad
\tau^{1}=\inf\{t\geq\gamma^{1} : x_{t} \not\in 2B\},
$$
$$
\gamma^{n+1}=\inf\{t\geq\tau^{n} : x_{t} \in \bar B\},\quad
\tau^{n+1}=\inf\{t\geq\gamma^{n+1} : x_{t} \not\in 2B\}.
$$
Then for any nonnegative Borel $f$ vanishing outside $B$
with $\|f\|_{L_{d}(B)}=1$
we have
$$ 
\int_{B}f   g_{1}(x)\,dx 
=E\int_{0}^{\infty}e^{-t}f ( x_{t})\,dt 
$$
$$
=\sum_{n=1}^{\infty}Ee^{- \gamma^{n}} 
E\Big(\int_{\gamma^{n}}^{\tau^{n}}e^{-(t-
 \gamma^{n} )}f( x_{t})\,dt
\mid \cF_{\gamma^{n}}\Big).
$$
Next we use the conditional  version of 
\eqref{10.4.7}
 to see that the conditional expectations
above are less than $NR $. After that
we use the conditional  version of Corollary
\ref{corollary 7.29.1} 
to get that
$$
R^{2}\leq N E\Big(\int_{\gamma^{n}}^{\tau^{n}}
e^{-(t-\gamma^{n})} \,dt
\mid \cF_{\gamma^{n}}\Big).
$$
Then we obtain
$$
\int_{B}f   g_{1}( x)\,dx 
\leq N  R ^{-1}
\sum_{n=1}^{\infty}E e^{-\gamma^{n}}
E\Big(\int_{\gamma^{n}}^{\tau^{n}}e^{-(t-\gamma^{n})} \,dt
\mid \cF_{\gamma^{n}}\Big)
$$
$$
=N  R ^{-1}
\sum_{n=1}^{\infty}E  
 \int_{\gamma^{n}}^{\tau^{n}}e^{- t } \,dt
$$
$$
\leq N  R ^{-1}
E \int_{0}^{\infty}e^{-t}I_{2B}( x_{t}) \,dt
=N  R ^{-1}\int_{2B}g_{1}( x)\,dx .
$$

The arbitrariness of $f$ 
implies that
$$
\Big(\dashint_{B}  g_{1} ^{d/(d-1)}(x)\,dx 
\Big)^{(d-1)/d}\leq N \dashint_{2B}  g_{1}(x)\,dx,
$$
and again it only remains to use Gehring's lemma
in case $p=d$.
For larger $p$ it suffices to use H\"older's inequality.
The theorem is proved.

Below by $d_{0}$ we mean the largest of $d_{0}$'s
from Theorems \ref{theorem 9.3.2} and \ref{theorem 10.4.3}
and first present improved parabolic Aleksandrov
estimates by following the
interpolation  arguments in Nazarov \cite{Na_15}.

\begin{lemma}
                                          \label{lemma 8.27.1}
Suppose that 
\begin{equation}
                                             \label{10.21.4}
p,q\in[1,\infty],\quad \frac{d_{0}}{p}+\frac{1}{q}=1.
\end{equation}
Then for any Borel $f(t,x)\geq0$ vanishing outside
a $Q\in\bQ$ of radius 1/2
\begin{equation}
                                             \label{8.27.1}
I:=E\int_{0}^{\infty}e^{-t}f(t,x_{t})\,dt
\leq N\|f\|_{L_{p,q}},
\end{equation}
where $N=N(\delta,d,\bar H,p ,p_{0},W)$.
\end{lemma}

Proof. 
If $p=d_{0}+1$, then $q=d_{0}+1$ and 
\eqref{8.27.1} follows from  \eqref{10.14.1}
and the fact that, obviously, $\|G_{1}\|_{L_{1,1}}=1$.
In other terms, for any Borel $f(t,x)\geq0$ vanishing outside
  $Q=(t_{0},t_{t})\times B_{R}(y) $  
$$
E\int_{0}^{\infty}e^{-t}f(t,x_{t})\,dt=\int_{Q}G_{1}(t,x)f(t,x)
\,dxdt
\leq N\|f\|_{L_{d_{0}+1,d_{0}+1}}.
$$

If $p=d_{0}$ and $q=\infty$ estimate \eqref{8.27.1}
follows from \eqref{10.4.6} since
$$
I\leq E\int_{0}^{\infty}e^{-t}\sup_{s\geq0}f(s,x_{t})\,dt
=\int_{\bR^{d}}g_{1}(x)\sup_{s\geq0}f(s,x)\,dx
$$
$$
\leq N\Big(\int_{\bR^{d}}\sup_{s\geq0}f^{d_{0}}(s,x)\,dx\Big)^{1/d_{0}}
=N\|f\|_{L_{d_{0},\infty}}.
$$

If $p=\infty$ and $q=1$
 $$
I\leq \int_{0}^{\infty} \sup_{x}f(t,x)\,dt=\|f\|_{L_{\infty,1}}.
$$

We will use these facts in an interpolation argument.
In case $\infty>p>d_{0}+1$ we have $p>q$ and
set $\beta=p/(d_{0}+1)$ and $\alpha=\beta/(\beta-1)$.
Take a nonnegative $g(t)$ such that $ \big(f(t,x)g(t)
\big)/g(t)=f(t,x)$ ($0/0=0$) and use H\"older's inequality to conclude
that $I\leq I_{1}I_{2}$, where
$$
I_{1}=\Big(\int_{0}^{\infty} g^{-\alpha}(t)
\,dt\Big)^{1/\alpha},
$$
$$
I_{2}=\Big(E\int_{0}^{\infty}e^{-t}g^{\beta}(t)f^{\beta}(t,x_{t})
\,dt\Big)^{1/\beta}
$$
$$
\leq N\Big(\int_{0}^{\infty}g^{(d_{0}+1)\beta}(t)
\Big(\int_{\bR^d}f^{(d_{0}+1)\beta}(t,x)\,dx\Big)\,dt
\Big)^{1/(d_{0}\beta+\beta)}.
$$
For $g$ found from
$$
  g^{-\alpha}(t)=g^{(d_{0}+1)\beta}(t)
\int_{\bR^d}f^{(d_{0}+1)\beta}(t,x)\,dx
$$
we get \eqref{8.27.1} and this takes care of the case
that $\infty>p>d_{0}+1$.

If $\infty>q>d_{0}+1$ we have $p<q$ and
set $\beta=q/(q-d_{0}-1)$ and $\alpha=\beta/(\beta-1)$.
Take a nonnegative $g(x)$ such that $ \big(f(t,x)g(x)
\big)/g(x)=f(t,x)$ ($0/0=0$) and use H\"older's 
inequality to conclude
that $I\leq I_{1}I_{2}$, where
$$
I_{1}=\Big(E\int_{0}^{\infty}e^{-t}g^{-\beta}(x_{t})\,dt
\Big)^{1/\beta}
\leq N\Big(\int_{\bR^{d}}g^{-d_{0}\beta}(x)\,dx
\Big)^{1/(d_{0}\beta)},
$$
$$
I_{2}=\Big(E\int_{0}^{\infty}e^{-t}g^{\alpha}(x_{t})
f^{\alpha}(t,x_{t})\,dt
\Big)^{1/\alpha}
$$
$$
\leq N\Big(\int_{\bR^{d}}g^{(d_{0}+1)\alpha}(x)
\Big(\int_{0}^{\infty}
f^{(d_{0}+1)\alpha}(t,x)\,dt\Big)dx
\Big)^{1/(\alpha d_{0}+\alpha )}.
$$
For $g$ found from
$$
g^{-d_{0}\beta}(x)=g^{(d_{0}+1)\alpha}(x)
 \int_{0}^{\infty}
f^{(d_{0}+1)\alpha}(t,x)\,dt 
$$
we get \eqref{8.27.1} and this proves the lemma.

Using this lemma instead of \eqref{9.29.2}
and just repeating the proof of Lemma \ref{lemma 8.22.1}
we come to a natural counterpart of the latter
and then by literally repeating the proof
of Theorem \ref{theorem 8.22.1} we come to the following  
result, that  is a version of Theorem 4.1 of Nazarov \cite{Na_15} in which
$d/p+1/q\leq 1$ that is stronger than ours,
but in which the assumption on $h$ is weaker.
A proper probabilistic 
version of Theorem 4.1 of Nazarov \cite{Na_15}
is found in \cite{Kr_20_2}. Set
$$
\bR^{d+1}_{+}=\{t\geq 0\}\cap\bR^{d+1}.
$$
\begin{theorem}
                                            \label{theorem 8.30.1}
Suppose   
\begin{equation}
                                                   \label{10.7.1}
p ,q \in [1,\infty],\quad
\nu:=1-\frac{d_{0}}{p }-\frac{1}{q }\geq 0.
\end{equation}
Then there is  $N=N(\delta,d,\bar H,p,q,p_{0},W)$ such that
for any $\lambda>0$ and   Borel nonnegative $f$  
we have 
\begin{equation}
                                                   \label{8.22.40}
 E\int_{0}^{\infty}e^{- \lambda t}
f(t,x_{t}) \,dt\leq 
N\lambda^{-\nu+(d-2d_{0})/(2p )}\|\Phi_{\lambda}^{1-\nu} f\|_
{L_{p ,q }(\bR^{d+1}_{+})},
\end{equation}
where $\Phi _{\lambda}(t,x)=\exp(- 
\sqrt\lambda (|x|+ \sqrt t)\theta/32)$.
In particular, if $f$ is independent of $t$, $p\geq d_{0}$,
 and $q=\infty$
$$
E\int_{0}^{\infty}e^{- \lambda t}
f( x_{t}) \,dt\leq 
N\lambda^{-1+d/(2p)}\|\bar \Phi_{\lambda}^{d_{0}/p} f\|_
{L_{p   }(\bR^{d } )},
$$
where $\bar \Phi _{\lambda}( x)=\exp(- 
\sqrt\lambda  |x| \theta/32)$.  
\end{theorem}

\begin{theorem}
                                          \label{theorem 9.7.1}
Assume that \eqref{10.7.1} holds.
Then
  for any
$n=1,2,...$, nonnegative Borel $f$ on $\bR^{d+1}_{+}$, and
 $T\in(0,\infty)$  we have
\begin{equation}
                                          \label{9.7.1}
E\Big[\int_{0}^{T}  
f(t,x_{t})\,dt\Big]^{n}\leq n!N^{n} 
T^{n\chi }\|\Phi^{(1-\nu)/n} _{1/T}
f\|^{n}_{L_{p,q}(\bR^{d+1}_{+}) },
\end{equation}
where $N=N(\delta,d,\bar H,p,q,p_{0},W)$ and 
$\chi=\nu+(2d_{0}-d)/(2p)$.  
 
\end{theorem}

Proof. The proof of this theorem proceeds by induction
on $n$ and is achieved by almost literally repeating
the proof of Theorem 2.6 of \cite{Kr_19_1}.
The induction hypothesis is that for all
$(t,x)\in \bR^{d+1}_{+}$ and $\kappa\in[0,1/n]$
\begin{equation}
                                                \label{10.28.1}
E\Big[\int_{0}^{T}  
f(t+s,x+x_{s})\,ds\Big]^{n}\leq n!N^{n} 
T^{n\chi }\Phi_{1/T}^{( \nu-1)\kappa n}
(t,x)\|\Phi^{( 1-\nu)\kappa} _{1/T}
f\|^{n}_{L_{p,q}(\bR^{d+1}_{+}) }.
\end{equation}

We will discuss in detail only the   case of $n=1$.
 In that case 
observe that, for any $\lambda>0$
$$
I:=E \int_{0}^{T}  
f(t+s,x+x_{s})\,ds\leq e^{\lambda T}
E \int_{0}^{\infty}e^{-\lambda s}  
f(t+s,x+x_{s})\,ds.
$$
It follows from Theorem \ref{theorem 8.30.1} 
 that for any   $\kappa\in[0,1]$
$$
  e^{-\lambda T}I \leq N \lambda^{-\chi}\|
f(t+\cdot,x+\cdot)
\Phi^{(1-\nu)\kappa} _{\lambda}\|_{L_{p,q}(\bR^{d+1}_{+}) }
$$
$$
\leq N \lambda^{-\chi}\Phi^{(\nu-1)\kappa}_{\lambda}(t,x)\|
f 
\Phi^{(1-\nu)\kappa} _{\lambda}\|_{L_{p,q}(\bR^{d+1}_{+}) },
$$
where the last inequality is due to the fact that,
for $(s,y)\in \bR^{d+1}_{+}$,
$\Phi  _{\lambda}(s,y)\leq \Phi  _{\lambda}(t+s,x+y)
\Phi ^{-1} _{\lambda}(t,x)$.
For $\lambda=1/T$ we get \eqref{10.28.1} with $n=1$. The theorem is 
proved.

Next theorem improves 
estimate \eqref{9.29.2} in what concerns
the restrictions on $p,q$.
\begin{theorem}
                                       \label{theorem 9.5.1}
Assume that \eqref{10.7.1} holds with $\nu=0$. 
 Then 
for any $R>0$ and Borel nonnegative $f$ 
given on $C_{R}$,
we have
\begin{equation}
                                                     \label{9.5.4}
E\int_{0}^{\tau_{R}}f(t,x_{t})\,dt\leq
NR^{(2d_{0}-d)/p}\|f\|_{L_{p,q}(C_{R})},
\end{equation}
where $N=N(\delta,d,\bar H,p, p_{0},W)$.

\end{theorem}

Proof.  
Self-similar transformations show that we may assume that $R=1$.
Also we may assume that $f$ is bounded and is zero outside $C_{1}$.
In that case denote by $\frM$ the set of stopping times
$\gamma\leq\tau:=\tau_{1}$, and set
$$
u_{\gamma}=E\Big[\int_{\gamma}^{\tau}f(t,x_{t})\,dt
\mid \cF_{\gamma}\Big],\quad  \bar u=\esssup_{\gamma\in\frM}u_{\gamma}.
$$
Observe that for any $\omega$  and $\lambda>0$
it holds that
$$
 \int_{\gamma}^{\tau}f(t,x_{t})\,dt=
\int_{\gamma}^{\tau}e^{-\lambda (t-\gamma)}f(t.x_{t})\,dt
$$
$$
+\lambda \int_{0}^{\infty}e^{-\lambda (t-\gamma)}
I_{\gamma\leq t<\tau}\Big[\int_{t}^{\tau}f(s,x_{s})\,ds\Big]\,dt.
$$
By the conditional  version of Theorem \ref{theorem 8.30.1} 
  (a.s.)
$$
E\Big[\int_{\gamma}^{\tau}e^{-\lambda (t-\gamma)}f(t,x_{t})\,dt
\mid \cF_{\gamma}\Big]
$$
$$
\leq
E\Big[\int_{\gamma}^{\infty}e^{-\lambda (t-\gamma)}f(t,x_{t})\,dt
\mid \cF_{\gamma}\Big]\leq 
N \lambda^{(d-2d_{0})/(2p )}\|f\|_{L_{p,q} }.
$$
Hence,
$$
u_{\gamma}\leq N \lambda^{(d-2d_{0})/(2p )}\|f\|_{L_{p,q} }
$$
$$
+\lambda E\Big[\int_{\gamma}^{\infty}e^{-\lambda (t-\gamma)}
I_{\gamma\leq t<\tau}E\Big[\int_{t}^{\tau}f(x_{s})\,ds
\mid \cF_{t}\Big]\,dt
\mid \cF_{\gamma}\Big],
$$
where the last term is dominated by
$$
\lambda \bar u E\Big[\int_{\gamma}^{\infty}e^{-\lambda (t-\gamma)}
I_{\gamma\leq t<\tau} \,dt
\mid \cF_{\gamma}\Big]\leq
\lambda \bar u E\Big[\int_{\gamma}^{\tau} \,dt
\mid \cF_{\gamma}\Big]\leq N_{1}\lambda \bar u
$$
(a.s.), where the last inequality follows from 
the  conditional  version of \eqref{9.29.2}.
Thus, (a.s.)
$$
u_{\gamma}\leq N \lambda^{(d-2d_{0})/(2p )}\|f\|_{L_{p,q} }
+N_{1}\lambda \bar u.
$$
Since $\gamma$ is arbitrary within $\frM$,
\begin{equation}
                                                       \label{9.5.7}
\bar u \leq N \lambda^{(d-2d_{0})/(2p )}\|f\|_{L_{p,q} }
+N_{1}\lambda \bar u
\end{equation}
(a.s.), and since $\bar u<\infty$ ($f$ is bounded), by taking $\lambda
=1/(2N_{1})$, we arrive at 
$$
\bar u\leq 
N  \|f\|_{L_{p,q} }.
$$
The theorem is proved.

\begin{remark}
                                            \label{remark 9.5.1}
Equation \eqref{9.5.7} implies
that  
$ d_{0}\geq d/2$. Of course, the example of the Wiener
process with no drift shows more than that, actually, $d_{0}> d/2$.
\end{remark}

Theorem \ref{theorem 9.5.1} allows us to prove It\^o's formula
for functions $u\in W^{1,2}_{p,q}(Q)$, where $Q$
is a domain in $\bR^{d+1}$ and 
$$
W^{1,2}_{p,q }(Q)=\{v: v, \partial_{t}v,
 Dv,  D^{2}v\in L_{p,q }(Q) \}
$$
with norm introduced in a natural way.
Before, the formula was known only for
(smooth, It\^o, and) $W^{1,2}_{d+1}$-functions
and processes with bounded drifts or
for $W^{2}_{d_{0}}$-functions in case the drift
of the process is dominated by $h(x_{t})$
with $h\in L_{d}$ (see \cite{Kr_19_1}).

The following extends Theorem 2.10.1  of \cite{Kr_77}.

\begin{theorem}
                                 \label{theorem 10.15.1}
Assume that \eqref{10.7.1} holds with $\nu=0$
and  $p<\infty$,
$q<\infty$. Let
$Q$ be a bounded domain in $\bR^{d+1}$, $0\in Q$,
$b$ be bounded, and  $u\in W^{1,2}_{p,q}(Q)\cap C(\bar Q)$. Then,
for $\tau$ defined as the first exit time of $(t,x_{t})$
from $Q$ and for all $t\geq0$,
$$
u(t\wedge\tau,x_{t\wedge\tau})
=u(0,0)+\int_{0}^{t\wedge\tau}D_{i}u(s,x_{s})\,dm^{i}_{s}
$$
\begin{equation}
                                      \label{10.15.1}
+\int_{0}^{t\wedge\tau}[
\partial_{t}u(s,x_{s})+ a^{ij}_{s}D_{ij}u(s,x_{s})
+b^{i}_{s}D_{i}u(s,x_{s})]\,ds
\end{equation}
and the stochastic integral above is a square-integrable
martingale.
\end{theorem}

Proof. First assume that $u$ is smooth and its derivatives
 are bounded. Then  
\eqref{10.15.1} holds by It\^o's formula and, moreover,
$$
E\int_{0}^{\tau}|Du(s,x_{s})|^{2}\,ds
\leq NE\Big(\int_{0}^{ \tau}D_{i}u(s,x_{s})\,dm^{i}_{s}\Big)^{2}
$$
$$
=NE\Big(u( \tau,x_{ \tau})-u(0,0)-
\int_{0}^{ \tau}[
\partial_{t}u(s,x_{s})+ a^{ij}_{s}D_{ij}u(s,x_{s})
+b^{i}_{s}D_{i}u(s,x_{s})]\,ds\Big)^{2}
$$
$$
\leq N\sup_{\bar Q}|u|+NE\Big(\int_{0}^{T}I_{Q}
\big(|\partial_{t}u|+|Du|+|D^{2}u|\big)(s,x_{s})\,ds\Big)^{2},
$$
where $T$ is the size of $Q$ in the $t$-direction.
In light of Theorem \ref{theorem 9.7.1} we conclude that
\begin{equation}
                                      \label{10.15.2}
E\int_{0}^{\tau}|Du(s,x_{s})|^{2}\,ds\leq N
\sup_{\bar Q}|u|+NT^{(2d_{0}-d)/p}\| \partial_{t}u,
 Du, D^{2}u\|_{L_{p,q}(Q)},
\end{equation}
where $N$ are independent of $u$ and $Q$.
Owing to Fatou's theorem,
this estimate is also true for those $u\in W^{1,2}_{p,q}(Q)
\cap C(\bar Q)$
that can be approximated uniformly and in the 
$W^{1,2}_{p,q}(Q)$-norm by smooth functions with bounded
derivatives (recall that $p<\infty$,
$q<\infty$). For our $u$
there is no guarantee that such approximation is possible.
However, mollifiers do such approximations
 in any subdomain $Q'\subset \bar Q'\subset
Q$. Hence, \eqref{10.15.2} holds for our $u$ if we replace $Q$
by $Q'$
(containing $(0,0)$). Setting $Q' \uparrow Q$  proves \eqref{10.15.2}
in the generals case and proves the last assertion
of the theorem.

After that \eqref{10.15.1} with $Q'$ in place of $Q$
is proved by routine approximation of $u$ by smooth 
functions. Setting $Q' \uparrow Q$ finally proves \eqref{10.15.1}.
The theorem is proved.

By using our Theorems \ref{theorem 8.20.1}
and \ref{theorem 9.5.1} (with $p=d+1$)
instead of Theorems 2.10 and 4.1 of \cite{Kr_19}
and arguing as in the proof of Theorem 4.4 of \cite{Kr_19}
we come to the following. Theorem \ref{theorem 11.8.1}
originated in \cite{KS_79} in case $b$ is bounded.
\begin{theorem}
                                     \label{theorem 11.8.1}
For any $\kappa\in(0,1)$ there is a 
function $q(\gamma)$, $\gamma\in(0,1)$,
depending only on $\kappa,\delta,d,\bar H, p_{0},W$, 
and, naturally, on $\gamma$,
such that for any $R\in(0,\infty)$, $x\in B_{\kappa R}$,
and closed $\Gamma\subset C_{ R}$ satisfying
$|\Gamma|\geq \gamma|C_{ R}|$ we have
$$
P (\tau_{\Gamma}(x) \leq  \tau_{   R}(x) )\geq q(\gamma),
$$
where $\tau_{\Gamma}(x) $ is the first time the process 
$(t, x+x_{t})$
hits $\Gamma$  and  
 $\tau_{ R}(x) $ is its first exit time from
$C_{ R}$. Furthermore, 
$q(\gamma)\to 1$ as $\gamma\uparrow 1$.
\end{theorem}

\mysection{The case of diffusion processes}
                                   \label{section 10.24.1}

Fix a constant $\delta\in(0,1)$ and by
$\bS_{\delta}$ denote the set of $d\times d$-symmetric
matrices whose eigenvalues are between
$\delta$ and $\delta^{-1}$. In this section
we impose the following.
\begin{assumption}
                                \label{assumption 10.12.1}

(i) On $\bR^{d+1}$
 we are given Borel measurable $\sigma(t,x)$
and $b(t,x)$ with values  in $\bS_{\delta}$
and in $\bR^{d}$ respectively.

(ii) We are given $p_{0},q_{0} \in[1,\infty)$  
satisfying \eqref{5.10.1} and a function $h(t,x)$
satisfying \eqref{8.19.1} for any $0<\tau\leq\rho^{2}$
and $(t,x)\in\bR^{d+1}$, where $\bar H$ is a fixed constant
and $W$ is a function satisfying $W(0+)=0  , W(1)=1$.

(iii) We have $|b|\leq h$.
\end{assumption}

Let $\Omega$ be the set of $\bR^{d+1}$-valued
 continuous function $(t_{0}+t,x_{t})$, $t_{0}\in \bR$,
defined for $t\in[0,\infty)$.
For $\omega=\{(t_{0}+t,x_{t}),t\geq0 \}$, define
$\sft_{t}(\omega)=t_{0}+t$, $x_{t}(\omega)=x_{t}$,
and set $\cN_{t}=\sigma((\sft_{s},x_{s}),s\leq t)$,
$\cN=\cN_{\infty}$.  
In the following theorem which is Theorem
6.1 of \cite{Kr_20_2} we use the terminology from
\cite{Dy_63}.

\begin{theorem}
                                            \label{theorem 4.27.1}
 On $\bR^{d+1}$ there exists a strong Markov process
$$
X=\{(\sft_{t},x_{t}),\infty,\cN_{t}, P_{t,x})
$$
such that the process 
$$
X_{1}=\{(\sft_{t},x_{t}),\infty,
\cN_{t+}, P_{t,x})
$$
 is Markov and for any $(t,x)\in\bR^{d+1}$
there exists a $d$-dimensional Wiener process $w_{t}$, $t\geq0$,
which is a Wiener process relative to $\bar \cN_{t}$,
where $\bar \cN_{t}$ is the completion of $\cN_{t}$
with respect to $P_{t,x}$, and such that with 
$P_{t,x}$-probability one, for
all $s\geq 0$, $\sft_{s}=t+s$ and 
\begin{equation}
                                               \label{10.16.3}
x_{s}=x +\int_{0}^{s}\sigma(t  +r,x_{r})\,dw_{r}
+\int_{0}^{s}b(t  +r,x_{r})\,dr.
\end{equation}
 \end{theorem}

\begin{remark}
                                      \label{remark 10.29.1}

 To be completely rigorous, to refer to \cite{Kr_20_2}
we should have $b\in L_{p_{0},q_{0}}$ (globally),
and   \eqref{8.19.1} is not needed.
But with our $b$, owing to estimates \eqref{10.28.2}
and \eqref{8.22.4}, the arguments in \cite{Kr_20_2}
only simplify and do not require $b\in L_{p_{0},q_{0}}$.
Still it is worth saying that the author believes that
under the conditions in \cite{Kr_20_2} Harnack's inequality
is true. Regarding the H\"older continuity of caloric
functions in the same setting we have no guesses.
The H\"older continuity seems to require some sort
of self-similarity and the $ L_{p_{0},q_{0}}$-norm
is not preserved under such transformations
if $p_{0},q_{0}$ are subject to \eqref{5.10.1}.
\end{remark}
 
\begin{theorem}
                                           \label{theorem 9.24.1}
For any $\lambda>0$, $ p,q $ satisfying 
\eqref{10.7.1}, and Borel nonnegative
$f(t,x)$  and for  
$$
 R_{\lambda}f(t,x):=E_{t,x}\int_{0}^{\infty}
e^{-\lambda s}f(t+s,x_{s})\,ds.
$$
we have
\begin{equation}
                                         \label{10.11.1}
\| R_{\lambda}f \|_{L_{p,q}(\bR^{d+1}_{+})}\leq N
\lambda^{-1} \| f \|_{L_{p,q}(\bR^{d+1}_{+})},
\end{equation}
where $N =N(\delta,d,\bar H,p,q,p_{0},W)$.

\end{theorem}

Proof. Self-similar transformations allow us to
concentrate on $\lambda=1$. 
By Theorem \ref{theorem 8.30.1}
(after appropriate shift of the origin in $\bR^{d+1}$)
we have
\begin{equation}
                                         \label{10.23.1}
R_{1}f (t,x)\leq 
N 
\|\Phi_{1}^{1-\nu} f (t+\cdot,x+\cdot)\|_{L_{p ,q }
(\bR^{d+1}_{+})}.
\end{equation}

Then, first assume that $p,q<\infty$. If $p\geq q$,
\eqref{10.23.1} implies that
\begin{equation}
                                         \label{11.3.1}
\int_{\bR^{d}}|R_{1}f(t,x)|^{p}\,dx
\leq N 
\int_{\bR^{d}}\Big(\int_{0}^{\infty}
F^{q/p}(t,s,x)\,ds\Big)^{p/q}dx ,
\end{equation}
where
$$
F(t,s,x)=\int_{\bR^{d}}\Phi^{(1-\nu)p}_{1}(s,y)f^{p} 
(t+s,x+y)\,dy.
$$
By Minkowski's inequality the   integral
on the right in \eqref{11.3.1} is dominated by
$$
\Big(\int_{0}^{\infty}\Big(\int_{\bR^{d}}F(t,s,x)\,dx\Big)^{q/p}
\,ds\Big)^{p/q},
$$
where
$$
\int_{\bR^{d}}F(t,s,x)\,dx=\int_{\bR^{d}}f^{p} (t+s,y)\,dy
\int_{\bR^{d}}\Phi^{(1-\nu)p}_{1}(s,y)\,dy
$$
$$
\leq N e^{-\mu \sqrt{ s}}
\int_{\bR^{d}}f^{p} (t+s,y)\,dy,
$$
with $\mu=\mu(\delta,p,q,\bar H,W) >0$. Below by $\mu$
we denote all such constants.
It follows that
$$
\int_{\bR^{d}}|R_{1}f(t,x)|^{p}\,dx
\leq N \Big(\int_{0}^{\infty}
e^{-\mu \sqrt{ s}}\Big(
\int_{\bR^{d}}f^{p}(t+s,y)\,dy\Big)^{q/p}ds\Big)^{p/q},
$$
$$
\|R_{1}f\|^{q}_{L_{p,q}(\bR^{d+1}_{+})}
\leq N \|f \|^{q}_{L_{p,q}(\bR^{d+1}_{+})}.
$$
 
If $q\geq p$,
$$
\int_{0}^{\infty}|R_{1}f(t,x)|^{q}\,dt
\leq N\int_{0}^{\infty}\Big(\int_{\bR^{d}}
 F^{p/q}\,(t,x,y) dy\Big)^{q/p}dt
$$
where
$$
F(t,x,y)=\int_{0}^{\infty}
\Phi^{(1-\nu)q}_{1}(s,y)f^{q} 
(t+s,x+y)\,ds.
$$
By Minkowski's inequality
$$
\Big(\int_{0}^{\infty}|R_{1}f(t,x)|^{q}\,dt\Big)^{p/q}
\leq N\int_{\bR^{d}}\Big(\int_{0}^{\infty}
F(t,x,y)\,dt\Big)^{p/q}dy,
$$
where
$$
\int_{0}^{\infty}
F(t,x,y)\,dt\leq
\int_{0}^{\infty}
 f^{q} 
( s,x+y)\,ds\int_{0}^{\infty}
\Phi^{(1-\nu)q}_{1}(s,y) \,ds
$$
$$
\leq Ne^{-\mu|y|}\int_{0}^{\infty}
 f^{q} 
( s,x+y)\,ds.
$$
Hence,
$$
\Big(\int_{0}^{\infty}|R_{1}f(t,x)|^{q}\,dt\Big)^{p/q}
\leq N\int_{\bR^{d}}e^{-\mu |y|}\Big(
\int_{0}^{\infty}
 f^{q} 
( s,x+y)\,ds\Big)^{p/q}\,dy,
$$
$$
\|R_{1}f\|_{L_{p,q}(\bR^{d+1}_{+})}^{p}
\leq N\|f \|_{L_{p,q}(\bR^{d+1}_{+})}^{p}
$$
and we again come to \eqref{10.11.1}.
This takes care of the case that $p,q<\infty$.

If $\nu=1$, then \eqref{10.11.1} follows from the fact that
$$
\sup_{\bR^{d+1}_{+}} R_{1}f\leq \sup_{\bR^{d+1}_{+}}  f.
$$
If $\nu<1$ but $p=\infty$ (and $q<\infty$),
then \eqref{10.11.1} follows from Minkowski's
inequality and the fact that
$$
\sup_{x}R_{1}f(t,x)\leq\int_{0}^{\infty}e^{-s}
\sup_{x}f(t+s,x)\,dt.
$$

If $\nu<1$ but $q=\infty$ (and $p<\infty$),
then \eqref{10.23.1} implies that
$$
\sup_{t\geq 0} (R_{1}f(t,x))^{p}\leq N \int_{\bR^{d}}
e^{-\mu |y|}\sup_{t\geq 0}f^{p}(t,x+y)\,dy.
$$
One sees that in this case \eqref{10.11.1}
follows after using Minkowski's
inequality once more.
The theorem is proved.

\mysection{An application to
first-order parabolic equations}
                                     \label{section 12.21.1}

We suppose that Assumption \ref{assumption 10.12.1}
is satisfied  and take 
$p,q\in[1,\infty)$ such that
$$
\frac{d}{p}+\frac{1}{q}=1,\quad p>q_{0}d.
$$
\begin{theorem}
                                  \label{theorem 12.20.2}
Let $Q$ be a bounded domain in $\bR^{d+1}$
with $0\in Q$ and let $u\in W^{1,2}_{p,q}(Q)
\cap C(\bar Q)$ be such that $u\leq 0$
on $\partial'Q$ and
$$
\partial_{t}u+b^{i}D_{i}u\leq 0
$$
in $Q$ (a.e.). Then $u(0)\leq0$.
\end{theorem}

This theorem is an immediate consequence of
the following result in which one need only
send $\varepsilon\downarrow 0$ and take into account
that $d/p+d^{2}q_{0}/(pp_{0})= q_{0}d/p<1$.

\begin{lemma}
                                     \label{lemma 12.20.2}
Under the assumptions of Theorem \ref{theorem 12.20.2}
for any $\varepsilon>0$ we have
\begin{equation}
                                       \label{12.20.3}
u(0)\leq N\varepsilon^{-d/p-d^{2}q_{0}/(pp_{0})}
\|\varepsilon\Delta u\|_{L_{p,q}(Q)},
\end{equation}
where $N$ is independent of $\varepsilon$.
\end{lemma}

Proof. Take $\nu>0$ and let $x_{t}$
be a solution of the equation
$$
x_{t}=\varepsilon w_{t}+\int_{0}^{t}
I_{(0,\nu)}(|b(s,x_{s})|)b(s,x_{s})\,ds
$$
on a probability space, where $w_{t}$ is a Wiener
process. By Theorem \ref{theorem 10.15.1}
(with $t$ larger than the diameter of $Q$)
 $$
u(0)\leq -E
\int_{0}^{ \tau}\big[
\partial_{t}u+ (1/2)\varepsilon^{2}
\Delta u 
+I_{(0,\nu)}(|b| )b^{i} D_{i}u \big]
(s,x_{s})\,ds.
$$
 
Here
$$
\partial_{t}u + (1/2)\varepsilon^{2}
\Delta u 
+I_{(0,\nu)}(|b| )b^{i} D_{i}u 
\geq I_{(0,\nu)}(|b| )\big(\partial_{t}u 
+b^{i} D_{i}u \big)
$$
$$
+ (1/2)\varepsilon^{2}\Delta u 
+I_{[\nu,\infty)}(|b| )\partial_{t}u 
\geq (1/2)\varepsilon^{2}\Delta u 
+I_{[\nu,\infty)}(|b| )\partial_{t}u .
$$
Also observe that with $I=(\delta^{ij})$
we have
$$
|b|\leq N(\det(\varepsilon^{2}I/2)^{1/p_{0}}
\varepsilon^{-2d/p_{0}}h,
$$
which by \eqref{9.29.2}
implies that
$$
E\int_{0}^{\tau}\varepsilon^{2d/p}\big|
(1/2)\varepsilon^{2}\Delta u 
+I_{[\nu,\infty)}(|b| )\partial_{t}u\big|(s,x_{s})\,ds
$$
$$
\leq N\big(1+\varepsilon^{-2d^{2}q_{0}/(pp_{0})}\big)
\|(1/2)\varepsilon^{2}\Delta u 
+I_{[\nu,\infty)}(|b| )\partial_{t}u\|_{L_{p,q}(Q)}.
$$

Hence,
$$
u(0)\leq N\varepsilon^{-2d/p}
 \big(1+\varepsilon^{-2d^{2}q_{0}/(pp_{0})}\big)
\|(1/2)\varepsilon^{2}\Delta u 
+I_{[\nu,\infty)}(|b| )\partial_{t}u\|_{L_{p,q}(Q)}.
$$
This is true with $N$ independent of $\varepsilon$
and $\nu$. By sending $\nu\to\infty$
we get \eqref{12.20.3} wit $\varepsilon^{2}$
in place of $\varepsilon$ which is irrelevant.
The lemma is proved.

The result of Theorem \ref{theorem 12.20.2}
is close to be sharp in the following sense.

\begin{example}
                                    \label{example 12.21.1}
Take $\varepsilon\in(0,1)$ and $p_{0},q_{0}\in[1,\infty)$
such that 
$d/p_{0}+1/q_{0}=1+\varepsilon$,
$p_{0}<q_{0}d$ (say $q_{0}\geq 2$),
 and $p_{0}\geq d$. Then it turns out that
there exists $b\in L_{p_{0},q_{0}}$, $p,q
\in[1,\infty)$ such that $d/p+1/q=1$ and $p>q_{0}d$,
a bounded domain $Q\subset\bR^{d+1}$ such that
$0\in Q$, and a function $u\in W^{1,2}_{p,q}(Q)
\cap C(\bar Q)$   such that $u\leq 0$
on $\partial'Q$ and
$$
\partial_{t}u+b^{i}D_{i}u\leq 0
$$
in $Q$ (a.e.), but $u(0)>0$.

To show this set
$$
\alpha=\frac{1-\varepsilon^{2}}{q_{0}},
\quad \beta=\frac{1-\varepsilon^{2}}{p_{0}}d
$$
  and observe that, since $p_{0}<q_{0}d$,
we have $\alpha<\beta$ and $1-\alpha>1-\beta$,
so that there exists $p$ satisfying
\begin{equation}
                                          \label{12.21.3}
\frac{d}{1-\beta}>p>\frac{d}{1-\alpha}.
\end{equation}
Here the left-hand side is strictly bigger than
$q_{0}d$ since 
$$
q_{0}(1-\beta)=q_{0}\Big(1-\frac{d}{p_{0}}+
\frac{\varepsilon^{2}}{p_{0}}d\Big)
=q_{0}\Big( \frac{1}{q_{0}}-\varepsilon+
\frac{\varepsilon^{2}}{p_{0}}d\Big)<1
$$
in light of $p_{0}\geq d$ and $\varepsilon\in(0,1)$.
Therefore, we can choose $p$ satisfying 
\eqref{12.21.3} and such that $p>q_{0}d$
as required. 
 
After that $q$ is also defined and we set
$$
u(t,x)=2-\exp\big(|t|^{1-\alpha}+|x|^{1+\beta}\big),
\quad Q=\{(t,x):u(t,x)>0\},
$$
$$
  b(t,x)=-\frac{1-\alpha}{1+\beta}\,\frac{1}{|t|^{\alpha}
|x|^{\beta}}\frac{x}{|x|}\,\text{sign}\,\,\! t.
$$
Since  $\alpha q_{0}<1$ and $\beta p_{0}<d$,
 we have $b\in L_{p_{0},q_{0}}(Q)$.
Also the inequality $\alpha q <1$,
guaranteeing that  $\partial_{t}u
\in L_{p_{0},q_{0}}(Q)$, is equivalent to 
the right inequality in \eqref{12.21.3},
whereas $p(1-\beta)<d$,
guaranteeing that  $D^{2}u
\in L_{p_{0},q_{0}}(Q)$, is equivalent to 
the left inequality in \eqref{12.21.3}.
Hence, $u\in W^{1,2}_{p_{0},q_{0}}(Q)$,
$u$ is also continuous, equals zero
on the whole boundary of $Q$, $u(0)=1$,
and, as is easy to see,
$\partial_{t}u+b^{i}D_{i}u=0$ apart from
the plane $t=0$.

\end{example}
 
{\bf Acknowledgment}. The author is sincerely
grateful to A.I. Nazarov for useful discussion
of the results.

\end{document}